\documentclass[oneside,american]{amsart}
\usepackage[T1]{fontenc}
\usepackage[latin9]{inputenc}
\usepackage{verbatim}
\usepackage{float}
\usepackage{amsthm}
\usepackage{amstext}
\usepackage{amssymb}
\usepackage{graphicx}

\makeatletter

\floatstyle{ruled}
\newfloat{algorithm}{tbp}{loa}
\providecommand{\algorithmname}{Algorithm}
\floatname{algorithm}{\protect\algorithmname}

\numberwithin{equation}{section}
\numberwithin{figure}{section}

\makeatother

\usepackage{babel}
\begin{document}

\title[Optimization via Separated Representations and the CTD]{Optimization via Separated Representations and the Canonical Tensor
Decomposition}

\author{Matthew J Reynolds$^{\star}$, Gregory Beylkin$^{\dagger}$, and
Alireza Doostan$^{\star}$}

\address{$^{\star}$ Department of Aerospace Engineering Sciences\\
429 UCB\\
University of Colorado at Boulder\\
Boulder, CO 80309\\
$^{\dagger}$ Department of Applied Mathematics\\
UCB 526\\
University of Colorado at Boulder\\
Boulder, CO 80309}

\thanks{This material is based upon work supported by the U.S. Department
of Energy Office of Science, Office of Advanced Scientific Computing
Research, under Award Number DE-SC0006402, and NSF grants DMS-1228359
and CMMI-1454601.}


\begin{abstract}
We introduce a new, quadratically convergent algorithm for finding
maximum absolute value entries of tensors represented in the canonical
format. The computational complexity of the algorithm is linear in
the dimension of the tensor. We show how to use this algorithm
to find global maxima of non-convex multivariate functions in separated form. 
We demonstrate the performance of the new algorithms on several
examples. 
\end{abstract}

\maketitle

\section{Introduction}

Finding global extrema of a multivariate function is an ubiquitous
task with many approaches developed to address this problem (see e.g. \cite{HOR-PAR:2013}).
Unfortunately, no existing optimization method can guarantee that
the results of optimization are true global extrema unless restrictive
assumptions are placed on the function. Assumptions on smoothness
of the function do not help since it is easy to construct an example
of a function with numerous local extrema ``hiding'' the location
of the true one. While convexity assumptions are helpful for finding global maxima, in practical applications
there are many non-convex functions. For non-convex functions various
randomized search strategies have been suggested and used but none
can assure that the results are true global extrema (see, e.g., \cite{SPALL:2005,HOR-PAR:2013}). 

We propose a new approach to the problem of finding global extrema
of a multivariate function under the assumption that the function
has certain structure, namely, a separated representation with a reasonably
small separation rank. While this assumption limits the complexity
of the function, i.e. number of independent degrees of freedom in
its representation, there is no restriction on its convexity. Furthermore, a large number of functions that do not appear to have
an obvious separated representation do in fact possess one, as was
observed in \cite{BEY-MOH:2002,BEY-MOH:2005,BE-GA-MO:2009}. In particular,
separated representations have been used recently to address curse
of dimensionality issues that occur when solving PDEs with random
data in the context of uncertainty quantification (see, e.g., \cite{CH-LA-CU:2011,DOO-IAC:2009,DO-IA-ET:2007,DO-VA-IA:2013,H-D-M-N:2014,NOUY:2007,NOUY:2008,NOUY:2010}).

We present a surprisingly simple algorithm that uses canonical tensor
decompositions (CTDs) for the optimization process. A canonical decomposition
of a tensor $\mathbf{U}\in\mathbb{R}^{M_{1}\times M_{2}\times\dots\times M_{d}}$
is of the form,
\[
\mathbf{U}=U\left(i_{1}\dots i_{d}\right)=\sum_{l=1}^{r}s_{l}u_{i_{1}}^{l}u_{i_{2}}^{l}\cdots u_{i_{d}}^{l},
\]
where $u_{i_{j}}^{l}$, $j=1,\dots,d$ are entries of vectors in $\mathbb{R}^{M_{d}}$
and $r$ is the separation rank. Using the CTD format to find entries
of a tensor with the maximum absolute value was first explored in
\cite{E-H-L-M-Z:2013} via an analogue of the matrix power method.
Unfortunately, the algorithm in \cite{E-H-L-M-Z:2013} has the same
weakness as the matrix power method: its convergence can be very slow,
thus limiting its applicability. Instead of using the power method,
we introduce a quadratically convergent algorithm based on straightforward
sequential squaring of tensor entries in order to find entries with
the maximum absolute value. 

For a tensor with even a moderate number of dimensions finding the
maximum absolute value of the entries appears to be an intractable
problem. A brute-force search of a $d$-directional tensor requires
$\mathcal{O}\left(N^{d}\right)$ operations, an impractical computational
cost. However, for tensors in CTD format the situation is different
since the curse of dimensionality can be avoided. 

Our approach to finding the entries with the maximum absolute value
(abbreviated as ``the maximum entry'' where it does not cause confusion)
is as follows: as in \cite{E-H-L-M-Z:2013}, we observe that replacing
the maximum entries with $1$'s and the rest of the entries with zeros
results in a tensor that has a low separation rank representation
(at least in the case where the number of such extrema is reasonably
small). To construct such a tensor, we simply square the entries via
the Hadamard product, normalize the result and repeat these steps
a sufficient number of times. The resulting algorithm is quadratically
convergent as we raise the tensor entries to the power $2^{n}$ after
$n$ steps. Since the nominal separation rank of the tensor grows
after each squaring, we use CTD rank reduction algorithms to keep
the separation rank at a reasonable level. The accuracy threshold
used in rank reduction algorithms limits the overall algorithm to
finding the global maximum within the introduced accuracy limitation.
However, in many practical applications, we can find the true extrema
as this accuracy threshold is user-controlled. These operations are
performed with a cost that depends linearly on the dimension $d$. 

The paper is organized as follows. In Section~\ref{sec:Background}
we briefly review separated representations of functions and demonstrate
how they give rise to CTDs. In Section~\ref{sec:A-quadratically-convergent}
we introduce the quadratically convergent algorithm for finding the
maximum entry of a tensor in CTD format and discuss the
selection of tensor norm for this algorithm. In Section~\ref{sec:CTD-based-global-opt}
we use the new algorithm to find the maximum absolute value of
continuous functions represented in separated form. In Section~\ref{sec:Numerical-examples}
we test both the CTD and separated representation optimization algorithms
on numerical examples. We provide our conclusions and a discussion
of these algorithms in Section~\ref{sec:Discussion-and-conclusions}.

\section{Background\label{sec:Background}}

\subsection{Separated representation of functions and the canonical tensor decomposition}

Separated representations of multivariate functions and operators
for computing in higher dimensions were introduced in \cite{BEY-MOH:2002,BEY-MOH:2005}.
The separated representation of a function $u(x_{1},x_{2},\dots,x_{d})$ is a natural extension
of separation of variables as we seek an approximation 
\begin{equation}
u(x_{1},\ldots,x_{d})=\sum_{l=1}^{r}s_{l}u_{1}^{(l)}(x_{1})\cdots u_{d}^{(l)}(x_{d})+\mathcal{O}\left(\epsilon\right),\label{eqn:sumsep}
\end{equation}
where $s_{l}>0$ are referred to as $s$-values. In this approximation
the functions $u_{j}^{(l)}(x_{j})$, $j=1,\dots,d$ are not fixed
in advance but are optimized in order to achieve the accuracy goal
with (ideally) a minimal \emph{separation rank} $r$. In (\ref{eqn:sumsep})
we set $x_{j}\in\mathbb{R}$ while noting that in general the variables
$x_{j}$ may be complex-valued or low dimensional vectors. Importantly,
a separated representation is not a projection onto a subspace, but
rather a nonlinear method to track a function in a high-dimensional
space using a small number of parameters. We note that the separation
rank indicates just the nominal number of terms in the representation
and is not necessarily minimal.

Any discretization of the univariate functions $u_{j}^{(l)}\left(x_{j}\right)$
in (\ref{eqn:sumsep}) with $u_{i_{j}}^{(l)}=u_{j}^{(l)}\left(x_{i_{j}}\right)$,
$i_{j}=1,\dots,M_{j}$ and $j=1,\dots,d$, leads to a $d$-dimensional
tensor $\mathbf{U}\in\mathbb{R}^{M_{1}\times\cdots\times M_{d}}$,
a canonical tensor decomposition (CTD) of separation rank $r_{u}$,
\begin{equation}
\mathbf{U}=U\left(i_{1},\dots,i_{d}\right)=\sum_{l=1}^{r_{u}}s_{l}^{u}\prod_{j=1}^{d}u_{i_{j}}^{(l)},\label{eq:sep-rep-algebraic}
\end{equation}
where the $s$-values $s_{l}^{u}$ are chosen so that each vector
$\mathbf{u}_{j}^{(l)}=\left\{ u_{i_{j}}^{(l)}\right\} _{i_{j}=1}^{M_{j}}$
has unit Euclidean norm $\Vert\mathbf{u}_{j}^{(l)}\Vert_{2}=1$ for
all $j,l$. The CTD has become one of the key tools in the emerging
field of numerical multilinear algebra (see, e.g., the reviews \cite{BRO:1997,TOM-BRO:2006,KOL-BAD:2009}).
Given CTDs of two tensors $\mathbf{U}$ and $\mathbf{V}$ of separation
ranks $r_{u}$ and $r_{v}$, their inner product is defined as 
\begin{equation}
\left\langle \mathbf{U},\mathbf{V}\right\rangle =\sum_{l=1}^{r_{u}}\sum_{l^{\prime}=1}^{r_{v}}s_{l}^{u}s_{l^{\prime}}^{v}\prod_{j=1}^{d}\left\langle \mathbf{u}_{j}^{(l)},\mathbf{v}_{j}^{(l^{\prime})}\right\rangle, \label{eq:tensor inner product}
\end{equation}
where the inner product $\left\langle \cdot,\cdot\right\rangle $
operating on vectors is the standard vector dot product. The Frobenius
norm is then defined as $\left\Vert \mathbf{U}\right\Vert _{F}=\sqrt{\left\langle \mathbf{U},\mathbf{U}\right\rangle }$.
Central to our optimization algorithm is the Hadamard, or point-wise,
product of two tensors represented in CTD format. We define this product
as,
\[
\mathbf{U}*\mathbf{V}=\left(U*V\right)\left(i_{1},\dots,i_{d}\right)=\sum_{l=1}^{r_{u}}\sum_{l^{\prime}=1}^{r_{v}}s_{l}^{u}s_{l^{\prime}}^{v}\prod_{j=1}^{d}\left(u_{i_{j}}^{(l)}\cdot v_{i_{j}}^{(l^{\prime})}\right).
\]

Common operations involving CTDs, e.g. sum or
multiplication, lead to CTDs with separation ranks that may be larger
than necessary for a user-specified accuracy. To keep computations with CTDs manageable, it is crucial to reduce
the separation rank. The separation rank
reduction operation for CTDs, here referred to as $\tau_{\epsilon}$,
is defined as follows: given a tensor $\mathbf{U}$ in CTD format
with separation rank $r{}_{u}$, 
\[
\mathbf{U}=U\left(i_{1},\dots,i_{d}\right)=\sum_{l=1}^{r_{u}}s_{l}^{u}\prod_{j=1}^{d}u_{i_{j}}^{(l)},
\]
and user-supplied error $\epsilon$, find a representation 
\[
\mathbf{V}=V\left(i_{1},\dots,i_{d}\right)=\sum_{l^{\prime}=1}^{r_{v}}s_{l^{\prime}}^{v}\prod_{j=1}^{d}v_{i_{j}}^{(l^{\prime})},
\]
with lower separation rank, $r_{v}<r_{u}$, such that $\left\Vert \mathbf{U}-\mathbf{V}\right\Vert <\epsilon\left\Vert \mathbf{U}\right\Vert $.
The workhorse algorithm for the separation rank reduction problem
is Alternating Least Squares (ALS) which was introduced originally
for data fitting as the PARAFAC (PARAllel FACtor) \cite{HARSHM:1970}
and the CANDECOMP \cite{CAR-CHA:1970} models. ALS has been used extensively
in data analysis of (mostly) three-way arrays (see e.g. the reviews
\cite{BRO:1997,TOM-BRO:2006,KOL-BAD:2009} and references therein).
For the experiments in this paper we use exclusively the randomized
interpolative CTD tensor decomposition, CTD-ID, described in \cite{BI-BE-BE:2015},
as a faster alternative to ALS.

\subsection{Power method for finding the maximum entry of a tensor\label{sub:The-power-method}}

The method for finding the maximum entry of a tensor in CTD format
in \cite{E-H-L-M-Z:2013} relies on a variation of the matrix power
method adapted to tensors. To see how finding entries with maximum
absolute value can be cast as an eigenvalue problem, we define the low
rank tensor $\mathbf{Y}$ with $1$s at the locations of such entries
in $\mathbf{U}$, and set $\lambda=\max_{i_{1},\dots,i_{d}}\left|U\left(i_{1},i_{2},\dots,i_{d}\right)\right|$.
The Hadamard product,
\begin{equation}
\mathbf{U}*\mathbf{Y}=\lambda\mathbf{Y}, \label{eq:eig problem}
\end{equation}
reveals the underlying eigenvalue problem. We describe the algorithm
of \cite{E-H-L-M-Z:2013} as Algorithm~\ref{alg:Power-method-algorithm}
and use it for performance comparisons. 
\begin{algorithm}[H]
\begin{raggedright}
To find the entries with maximum absolute value of a tensor $\mathbf{U}\in\mathbb{R}^{N_{1}\times N_{2}\times\dots\times N_{d}}$
in CTD format, given a tolerance $\epsilon>0$ for the reduction algorithm
$\tau_{\epsilon}$, we define\\
~\\

\par\end{raggedright}

\begin{raggedright}
$\mathbf{Y}_{0}=\left(\prod_{i=1}^{d}\frac{1}{N_{i}}\right)\left(\mathbf{1}_{1}\circ\mathbf{1}_{2}\circ\dots\circ\mathbf{1}_{d}\right)$,
where $\mathbf{1}_{i}=\left(1,1,\dots,1\right)^{T}\in\mathbb{R}^{N_{i}}$, 
\par\end{raggedright}

\begin{raggedright}
and iterate (up to $k_{\max}$ times) as follows:\\
~\\
\par\end{raggedright}

\begin{raggedright}
$\mathbf{for}\,\, k=1,2,\dots,k_{\max}$ $\mathbf{do}$
\par\end{raggedright}
\begin{enumerate}
\item \begin{raggedright}
$\quad\mathbf{Q}_{k}=\mathbf{U}*\mathbf{Y}_{k-1},\,\,\lambda_{k}=\left\langle \mathbf{Y}_{k-1},\mathbf{Q}_{k}\right\rangle ,\,\,\mathbf{Z}_{k}=\mathbf{Q}_{k}/\left\Vert \mathbf{Q}_{k}\right\Vert _{F}$ \\[-.15in]
\par\end{raggedright}
\item \begin{raggedright}
$\quad\mathbf{Y}_{k}=\tau_{\epsilon}\left(\mathbf{Z}_{k}\right)$
\par\end{raggedright}
\end{enumerate}
\begin{raggedright}
$\mathbf{end}\,\,\mathbf{for}$
\par\end{raggedright}

\protect\caption{\label{alg:Power-method-algorithm} }

\end{algorithm}
\noindent The resulting tensor
$\mathbf{Y}_{k}$ ideally has a low separation rank with zero entries
except for a few non-zeros indicating the location of the extrema.
In the case of a single extremum the tensor $\mathbf{Y}_{k}$ has
separation rank $1$.

\section{A quadratically convergent algorithm for finding maximum entries
of tensors\label{sec:A-quadratically-convergent}}

In order to construct  $\mathbf{Y}$ from (\ref{eq:eig problem}),
we simply square all entries
of the tensor, normalize the result, and repeat these steps a number
of times. Fortunately, the CTD format allows us to square tensor entries
using the Hadamard product. These operations increase the separation
rank, so the squaring step is (usually) followed by the application
of a rank reduction algorithm such as ALS, CTD-ID, or their combination.
\begin{algorithm}[H]
\begin{raggedright}
To find the entries with maximum absolute value of a tensor $\mathbf{U}\in\mathbb{R}^{N_{1}\times N_{2}\times\dots\times N_{d}}$
in CTD format, given a reduction tolerance $\epsilon>0$ for the reduction
algorithm $\tau_{\epsilon}$, we define\\
~\\

\par\end{raggedright}

\begin{raggedright}
$\mathbf{Y}_{0}=\mathbf{U}/\left\Vert \mathbf{U}\right\Vert _{F}$,
and iterate (up to $k_{\max}$ times) as follows:\\
~\\

\par\end{raggedright}

\begin{raggedright}
$\mathbf{for}\,\, k=1,2,\dots,k_{\max}$ $\mathbf{do}$
\par\end{raggedright}
\begin{enumerate}
\item \begin{raggedright}
$\quad\mathbf{Q}_{k}=\mathbf{Y}_{k-1}*\mathbf{Y}_{k-1}$\\[-.15in]
\par\end{raggedright}
\item \begin{raggedright}
$\quad\mathbf{Y}_{k}=\tau_{\epsilon}\left(\mathbf{Q}_{k}\right)$,
$\mathbf{Y}_{k}=\mathbf{Y}_{k}/\left\Vert \mathbf{Y}_{k}\right\Vert _{F}$ 
\par\end{raggedright}
\end{enumerate}
\begin{raggedright}
$\mathbf{end}\,\,\mathbf{for}$\\

\par\end{raggedright}

\protect\caption{\label{alg:squaring_algorithm} }
\end{algorithm}
 If we identify the largest (by absolute value) entry as $a=\left|F\left(j_{1},j_{2},\dots,j_{d}\right)\right|$
and the next largest as $b$, the ratio $b/a<1$ decreases rapidly,
\[
\left(\frac{b}{a}\right)^{2^{j}}\le\epsilon
\]
and the total number of iterations, $j$, to reach $\epsilon$ is
logarithmic, 
\[
j\le\log_{2}\left(\frac{\log\left(\epsilon\right)}{\log\left(b\right)-\log\left(a\right)}\right).
\]
Since this ratio is expected to be reasonably small for most entries,
the tensor will rapidly become sparse after only a small number of
iterations $j$. The quadratic rate of convergence of Algorithm~\ref{alg:squaring_algorithm}
is the key difference with the algorithm in \cite{E-H-L-M-Z:2013}
whose convergence rate is at best linear.

\subsection{Termination conditions for Algorithm~\ref{alg:squaring_algorithm}}

We have explored three options to terminate the iteration in Algorithm~\ref{alg:squaring_algorithm}.
First, the simplest, is to fix the number of iterations in advance
provided some knowledge of the gap between the entry with the maximum
absolute value and the rest of the entries is available.
The second option is to define $\lambda_{k}=\left\langle \mathbf{Y}_{k},\mathbf{U}\right\rangle $
and terminate the iteration once the rate of its decrease becomes
small, i.e., $\left(\lambda_{k}-\lambda_{k-1}\right)/\lambda_{k-1}<\delta$,
where $\delta>0$ is a user-supplied parameter. The third option is
to observe the separation rank of $\mathbf{Y}_{k}$ and terminate
the iteration when it becomes smaller than a fixed, user-selected
value. In Section~\ref{sec:Numerical-examples} we indicate which
termination condition is used in our numerical examples. 

\subsection{Selection of tensor norm\label{sub:Selection-of-tensor-norm}}

Working with tensors, we have a rather limited selection of norms
that are computable. The usual norm chosen for work with tensors is
the Frobenius norm which may not be appropriate in all situations
where the goal is to find the maximum absolute value entries. Unfortunately,
the Frobenius norm is only weakly sensitive to changes
of individual entries. The alternative $s$-norm (see \cite[Section 4]{BI-BE-BE:2015}),
i.e., the largest $s$-value of the rank one separated approximation
to the tensor in CTD format, is better in some situations. In particular,
it allows the user to lower the tolerance to below single precision
(within a double precision environment). The improved control over
truncation accuracy makes this norm more sensitive to individual entries
of the tensor. We note that while it is straightforward to compute
the $s$-value of a rank one separated approximation via a convergent
iteration, theoretically it is possible that this computed number
is not the largest $s$-value. However, we have not encountered such
a situation using the $s$-norm in practical problems.

\subsection{Numerical demonstration of convergence}

We construct an experiment using random tensors, with low separation
rank ($r=4$ in what follows), in dimension $d=6$ with $M=32$ samples
in each direction. The input CTD is formed by adding a rank three
CTD and a rank one CTD together. The rank three CTD is formed using
random factors ($\mathbf{u}_{j}^{(l)}$ in (\ref{eq:sep-rep-algebraic}))
whose samples were drawn from a uniform distribution on $\left[.9,1\right]$
(the $s$-values were set to $1$). The rank one CTD has all zeros
except for a spike added in at a random location. The size of this
spike is chosen to make the maximum entry have a magnitude of $3.5$.
For the experiments in this section we set the reduction tolerance
of $\tau_{\epsilon}$ in Algorithm~\ref{alg:squaring_algorithm}
to $\epsilon=10^{-6}$ and use the Frobenius norm.

\begin{figure}[ht]
\begin{centering}
\includegraphics{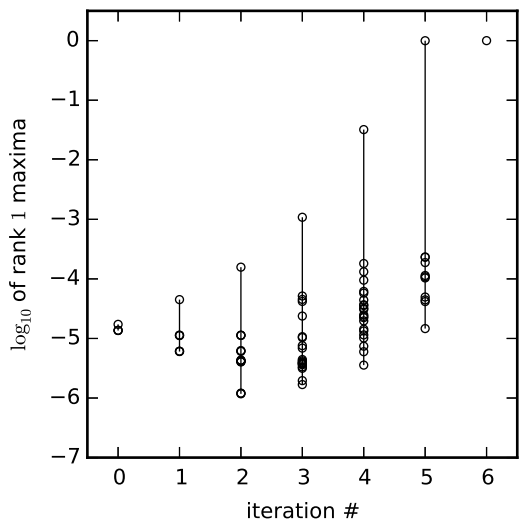}

\par\end{centering}

\protect\caption{\label{fig:factor_maxima_comparison}Comparing the maximum absolute
values of entries from rank one terms of a CTD. After $6$ iterations
 Algorithm~\ref{alg:squaring_algorithm} converges to a rank $1$
tensor.}
\end{figure}

We observe that as the iteration proceeds, the maximum entry of one
of the CTD terms eventually becomes dominant relative to the rest.
The iteration continues until the rest of the terms are small enough
to be eliminated by the separation rank reduction algorithm $\tau_{\epsilon}$.
This is illustrated in Figure~\ref{fig:factor_maxima_comparison}
where we display the maximum entry from all rank one terms in the
CTD for each iteration. Before the first iteration (iteration $0$
in Figure~\ref{fig:factor_maxima_comparison}), the gaps between
the maxima of the rank one terms are not large. As we iterate, the
maximum entry of one of the rank one terms rapidly separates itself
from the others and, by iteration $6$, the reduction algorithm removes
all terms except one. This remaining term has an entry $1$ corresponding
to the location of the maximum entry of the original tensor and all
other entries are $0$.

\begin{figure}[ht]
\begin{centering}
\includegraphics{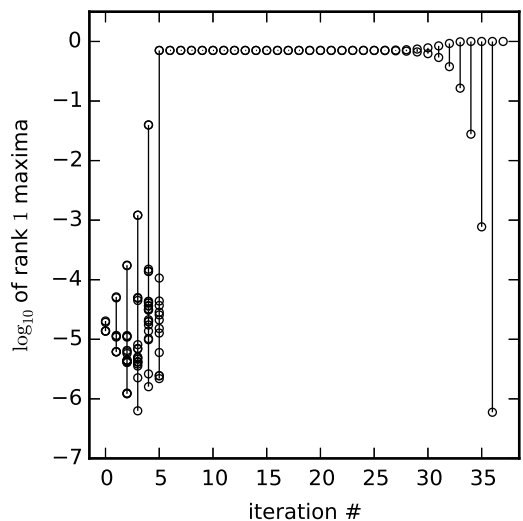}

\par\end{centering}

\protect\caption{\label{fig:factor_maxima_comparison-mult}In this example there are
two maxima of the same size. The values at these locations converge
to $1/\sqrt{2}$ and the rest of the entries converge to $0$. If allowed
to run long enough, the value at one of the locations becomes dominant
due to the truncation error of the rank reduction algorithm. The other
location is then eliminated by the algorithm and only one term remains.}
\end{figure}

Another interesting case occurs when a tensor has multiple maximum
entry candidates, either close to or exactly the same in magnitude.
To explore this case, we construct an experiment similar to the previous
one, namely, find the maximum absolute value of a CTD of dimension
$d=6$ and $M=32$ samples in each direction. This CTD is constructed
with the same rank $3$ random CTD as above, but this time we added
in two spikes at random locations so that there are two maximal spikes
of absolute value $3.5$. In cases such as this, Algorithm~\ref{alg:squaring_algorithm}
can be run for a predetermined number of iterations, after which the
factors are examined and a subset of candidates can be identified
and explicitly verified. For example, in Figure~\ref{fig:factor_maxima_comparison-mult},
we observe that by iteration $6$ there are only two candidates for
the location of the maximum. 
As a warning we note that allowing the code to continue to run many
more iterations will leave only one term due to the finite
reduction tolerance $\epsilon$ of $\tau_{\epsilon}$ in Algorithm~\ref{alg:squaring_algorithm}.
Indeed, in Figure~\ref{fig:factor_maxima_comparison-mult} we observe
the maximum entries of the two rank one terms begin to split around
iteration number $26$, and by iteration $37$ only one term remains.
This example illustrates a potential danger of insisting on a single
term termination condition. 

\section{CTD-based global optimization of multivariate functions\label{sec:CTD-based-global-opt}}

We describe an approach using Algorithm~\ref{alg:squaring_algorithm} for the global optimization of functions that admit low rank separated
representations.
We consider two types of problems: first, where the
objective function is already in a separated form, ready to yield a
CTD appropriate for use with Algorithm~\ref{alg:squaring_algorithm}.
The second, more general problem, is to construct a separated representation
given data or an analytic expression of the function, from which a
CTD can then be built. 

In the first problem the objective function is already represented
in separated form as in (\ref{eqn:sumsep}), and an interpolation
scheme is associated with each direction. Such problems have been
subject to growing interest due to the use of CTDs for solving operator
equations, e.g., deterministic or stochastic PDE/ODE systems (see,
for example, \cite{BEY-MOH:2002,BEY-MOH:2005,DOO-IAC:2009,NOUY:2010,CH-LA-CU:2011,CH-VE-KA:2016}). 

In the second problem the goal is to optimize a multivariate function
where no separated representation is readily available. What typically
is available is a data set of function values (scattered or on a grid).
In such cases our recourse is to first construct a separated representation of the
underlying multivariate objective function and then re-sample it in
order to obtain a CTD of the form (\ref{eq:sep-rep-algebraic}). The
construction of a separated representation given a set of function
values can be thought of as a regression problem, see \cite{BE-GA-MO:2009,DO-VA-IA:2013}.
The algorithms in these papers use function samples
to build a separated representation of form (\ref{eqn:sumsep}). Following
\cite{BE-GA-MO:2009}, an interpolation scheme is set up in each direction
and the ALS algorithm is used to reduce the problem to regression
in a single variable. The univariate interpolation scheme may use
a variety of possible bases, as well as nonlinear approximation techniques.
Once a separated representation is constructed, the CTD for finding
the global maxima is then obtained to satisfy the local interpolation
requirement stated below. The resulting CTD is then used as input into Algorithm~\ref{alg:squaring_algorithm}
to produce an approximation of the global maxima.

If the objective function is given analytically, it is preferable
to use analytic techniques to approximate the original function via
a separated representation. An example of this approach is included
in Section~\ref{sub:Optimization-example} below.

We note that sampling becomes an important issue when finding global
maxima of a function via Algorithm~\ref{alg:squaring_algorithm}.
The key sampling requirement is to ensure that the objective function
can be interpolated at an arbitrary point from its sampled values
up to a desired accuracy. This implies that sampling rate must
depend on the local behavior of the function, e.g., the sampling rate
is greater near singularities. For example, if a function has algebraic
singularities, an efficient method to interpolate is to use wavelet
bases. In such a case we emphasize the importance of using interpolating
scaling functions, i.e., scaling functions whose coefficients are
function values or are simply related to those. Since the global maxima
are typically searched for in finite domains, the multiresolution
basis should also work well in an interval (or a box in higher dimensions).
This leads us to suggest the use of Lagrange interpolating polynomials
within an adaptive spatial subdivision framework. We note that rescaled
Lagrange interpolating polynomials form an orthonormal basis on an
interval and the spatial subdivision is available as multiwavelet
bases (see e.g., \cite{ALPERT:1993,A-B-G-V:2002}).

Finding the maximum entries of a CTD yields only an approximation
of the maxima and their locations for the function. However, due to
the interpolation requirements stated above these are high-quality
approximations so that, if needed, a local optimization algorithm
can then be used to find the true global maxima. Since in this case
we use local optimization, it is preferable to oversample the CTD
representation of a function, i.e. to use large $M_{j}$'s in (\ref{eq:sep-rep-algebraic}).
This will not result in a prohibitive cost as the rank reduction algorithms,
e.g. ALS, that operate on CTDs scale linearly with respect to the
number of samples in each direction \cite{BEY-MOH:2005}. However,
the total number of samples needed to satisfy the local interpolation
requirement may be large and, thus, additional steps to accelerate
the reduction algorithm $\tau_{\epsilon}$ in Algorithm~\ref{alg:squaring_algorithm}
may be required. One such technique consists of using the $QR$
factorization as explained in \cite[Section 2.4]{BI-BE-BE:2015}.

\section{Numerical examples\label{sec:Numerical-examples}}

\subsection{Comparison with power method algorithm from Section~\ref{sub:The-power-method} }

To test Algorithm~\ref{alg:squaring_algorithm}, we construct tensors
from random factors and find the entries with the largest absolute
value. For these examples we select dimension $d=8$ and set $M=32$
samples in each direction. The test CTDs are constructed by adding
rank four and rank one CTDs together. The rank four CTDs are composed
of random factors (vectors $\mathbf{u}_{j}^{(l)}$ in (\ref{eq:sep-rep-algebraic}))
whose entries are drawn from the uniform distribution on $\left[.9,1\right]$
and the corresponding $s$-values are set to $1$. The rank one CTDs
contain all zeros except for a magnitude $4$ spike added in at a
random location. The CTDs in these examples have a unique maximum
entry, and we chose to terminate the algorithm when the output CTD
reaches separation rank one. We run tests for both Algorithm~\ref{alg:squaring_algorithm}
and Algorithm~\ref{alg:Power-method-algorithm}, using for the rank
reduction operation $\tau_{\epsilon}$ the CTD-ID algorithm with accuracy
threshold $\epsilon=10^{-6}$ in the $s$-norm.
%
\begin{figure}[ht]
\begin{centering}
\includegraphics{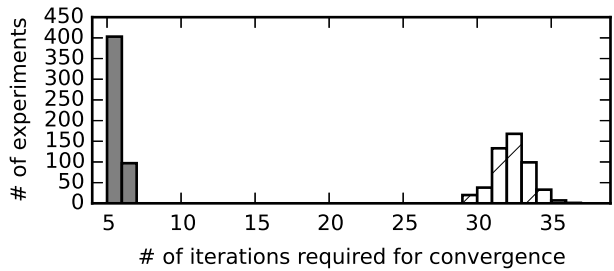}
\par\end{centering}

\protect\caption{\label{fig:random_iter_plots}Histograms illustrating the number of
iterations required for convergence to a rank one CTD by the power
method optimization (hatched pattern) and the squaring method optimization
(solid gray).}
\end{figure}

\begin{figure}[ht]
\begin{centering}
\includegraphics{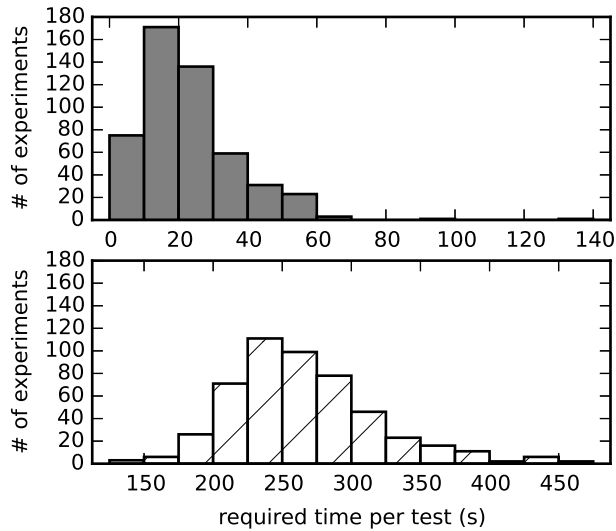}
\par\end{centering}

\protect\caption{\label{fig:random_time_plots}Histograms illustrating the computation
times, in seconds, required for convergence to a rank one CTD by the
power method optimization (hatched pattern) and the squaring method
optimization (solid gray).}
\end{figure}

A histogram displaying the required number of iterations to converge
for $500$ individual tests is shown in Figure~\ref{fig:random_iter_plots}.
While both algorithms find the correct maximum location in all tests,
Algorithm~\ref{alg:squaring_algorithm} outperforms Algorithm~\ref{alg:Power-method-algorithm}
in terms of the number of iterations required. While Algorithm~\ref{alg:squaring_algorithm}
requires much fewer iterations than the power method for the same
problem, a lingering concern may be that the squaring process in Algorithm~\ref{alg:squaring_algorithm}
requires the reduction of CTDs with much larger separation ranks.
In Figure~\ref{fig:random_time_plots} we show histograms of computation
times. We observe that despite reducing larger separation rank CTDs,
the computation times using Algorithm~\ref{alg:squaring_algorithm}
are significantly smaller. These computation times are for our
MATLAB code running on individual cores of a Intel Xeon $2.40$ GHz
processor.

\subsection{\label{sub:Optimization-example}Optimization example}

We reiterate that separated representations of general functions can 
be obtained numerically using only function evaluations, as described
in \cite{BE-GA-MO:2009,DO-VA-IA:2013}. Once such a representation
is obtained, an appropriate sampling yields a tensor in CTD format. 
In some cases, as in this example, this can be accomplished analytically.

As an example of constructing separated representations for optimization,
we consider Ackley's test function \cite{BACK:1996}, commonly used
for testing global optimization algorithms,
\begin{equation}
u\left(\mathbf{x}\right)=a\,e^{-b\left(\frac{1}{d}\sum_{i=1}^{d}x_{i}^{2}\right)^{1/2}}+e^{\frac{1}{d}\sum_{i=1}^{d}\cos\left(cx_{i}\right)},\label{eq:ackley_fun}
\end{equation}
where $d$ is the number of dimensions and $a$, $b$, and $c$ are
parameters. For our tests we set $d=10$, $a=20$, $b=0.2$, and $c=2\pi$.
We choose this test function for two reasons: first, it has many local
maxima, thus a local method without a good initial guess may not find the correct maximum. Second, one of the terms
is not in separated form, hence we must construct a separated approximation
of the function first. 

The true maximum of (\ref{eq:ackley_fun}) occurs at $\mathbf{x}=\mathbf{0}$
with a value of $u\left(\mathbf{0}\right)=a+e$. The first term in
(\ref{eq:ackley_fun}) is radially symmetric which we approximate
with a linear combination of Gaussians using the approximation of
$e^{-xy}$ from \cite[Section 2.2]{BEY-MON:2010},
\begin{equation}
G_{e}\left(x\right)=\frac{hb}{2\sqrt{\pi d}}\sum_{j=0}^{R}\exp\left(-\frac{b^{2}}{4d}e^{s_{j}}-x^{2}e^{s_{j}}+\frac{1}{2}s_{j}\right),\label{eq:expansion-by-gaussians}
\end{equation}
where $s_{j}=s_{start}+jh$, and $h,\,\, s_{start},$ and $R$ are
parameters chosen such that given $\epsilon,\delta>0$,
\begin{equation}
\left|G_{e}\left(x\right)-e^{-\frac{b}{\sqrt{d}}x}\right|\le\epsilon\label{eq:gauss-exp-difference}
\end{equation}
for $0<\delta<x<\infty$. We arrive at the approximation of the first
term in (\ref{eq:ackley_fun}),
\begin{equation}
\left|e^{ -b \left(\frac{1}{d}\sum_{i=1}^{d}x_{i}^{2}\right)^{1/2}}-\sum_{j=0}^{R}w_{j}\prod_{i=1}^{d}\exp\left(-x_{i}^{2}e^{s_{j}}\right)\right|\le\epsilon\label{eq:SR-of-radial-exponential}
\end{equation}
where $w_{j}=\frac{hb}{2\sqrt{\pi d}}\exp\left(-\frac{b^{2}}{4d}e^{s_{j}}+\frac{1}{2}s_{j}\right)$.

To correctly sample (\ref{eq:ackley_fun}), we first consider two
terms of the function separately, the radially-symmetric exponential
term, and the exponential cosine term which is already in separated form, $\exp\left(\frac{1}{d}\sum_{i=1}^{d}\cos\left(cx_{i}\right)\right)=\prod_{i=1}^{d}\exp\left(\frac{1}{d}\cos\left(cx_{i}\right)\right)$.
These terms require different sampling strategies. The sampling rate
of the radially-symmetric exponential term near the origin should
be significantly higher than away from it. On the other hand, uniform sampling is sufficient for the
exponential cosine term. We briefly discuss
the combined sampling of these terms to satisfy both requirements.

To sample the radially-symmetric exponential term in (\ref{eq:ackley_fun}),
we use the expansion by Gaussians (\ref{eq:expansion-by-gaussians})
as a guide. We choose $10$ equally-spaced points centered at zero
and sample the Gaussian $e^{-x^{2}}$. For each Gaussian in (\ref{eq:expansion-by-gaussians}),
we scale these points by the factor $e^{s_{j}/2}$ in order to consistently
sample each Gaussian independent of its scale. Since (\ref{eq:expansion-by-gaussians})
includes a large number of rapidly decaying Gaussians, too many of
the resulting samples will concentrate near the origin. Therefore,
we select a subset as follows: first, we keep all samples from the
sharpest Gaussian, $x_{0l}$, $l=1,\dots,10$. We then take the samples
from second sharpest Gaussian and keep only those outside the range
the previous samples, i.e., we keep $x_{1m}$ such that $\left|x_{1m}\right|>\max_{l}\left|x_{0l}\right|$,
$m=1,\dots,10$. Repeating this process for all terms $j=0,1,\dots,R$
in expansion (\ref{eq:expansion-by-gaussians}) yields samples for
the radially-symmetric exponential term in (\ref{eq:ackley_fun}).

We oversample the exponential cosine term using $16$ samples per
oscillation. Notice that this sampling strategy is not sufficient
for the radially-symmetric exponential term in (\ref{eq:ackley_fun})
(and vice-versa). To combine these samples, we first find the radius
where the sampling rate of the radially-symmetric exponential term
in (\ref{eq:ackley_fun}) drops below the sampling rate chosen to
sample the exponential cosine term. We then remove all samples whose
coordinates are greater in absolute value than the radius, and replace
them with samples at the chosen rate of $16$ samples per oscillation.
The resulting set of samples consists of $M=1080$ points for each
direction $j=1,\dots,10$.

Using expansion (\ref{eq:expansion-by-gaussians}) leads to a $77$
term CTD representation for the chosen accuracy $\epsilon=10^{-8}$
and parameter $\delta=3\cdot10^{-6}$ in (\ref{eq:gauss-exp-difference}).
An initial reduction of this CTD using the CTD-ID algorithm with accuracy
threshold $\epsilon=10^{-6}$ in the Frobenius norm produces a CTD
with separation rank $11$ representing the function $u\left(\mathbf{x}\right)$.
We use this CTD as an input into Algorithm~\ref{alg:squaring_algorithm},
yielding a rank one solution after $24$ iterations. We set the reduction
tolerance in Algorithm~\ref{alg:squaring_algorithm} for this experiment
to $\epsilon=10^{-6}$ and find the maximum entry of the tensor located
at a distance $2.13\times10^{-3}$ away from the function's true maximum
location, the origin. Using the output as the initial guess in a local,
gradient-free, optimization algorithm (in our case a compass search,
see e.g. \cite{KO-LE-TO:2003}) yields the maximum value with relative
error of $1.03\times10^{-7}$ located at a distance $1.84\times10^{-6}$
away from the origin.

\section{Discussion and conclusions\label{sec:Discussion-and-conclusions}}
Entries of a tensor in CTD format with the largest absolute value
can be reliably identified by the new quadratically convergent algorithm,
Algorithm~\ref{alg:squaring_algorithm}. This algorithm, in turn,
can be used for global optimization of multivariate functions that
admit separated representations. Unlike algorithms based on random
search strategies (see e.g. \cite{HOR-PAR:2013}), our approach allows
the user to estimate the potential error of the result. The computational
cost of Algorithm~\ref{alg:squaring_algorithm} is dominated by the
cost of reducing the separation ranks of intermediate CTDs, not the
dimension of the optimization problem. If ALS is used to reduce the
separation rank and $d$ is the dimension, $r$ is the maximum separation
rank (after reductions) during the iteration, and $M$ is the maximum
number of components in each direction, then the computational cost
of each rank reduction step can be estimated as $\mathcal{O}\left(r^{4}\cdot M\cdot d\cdot N_{iter}\right)$,
where $N_{iter}$ is the number of iterations required by the ALS algorithm to converge. The computational cost
of the CTD-ID algorithm is estimated as $\mathcal{O}\left(r^{3}\cdot M\cdot d\right)$
and, if it is used instead of ALS, the reduction step is faster by
a factor of $\mathcal{O}\left(r\cdot N_{iter}\right)$ \cite{BI-BE-BE:2015}.
We note that, while linear in dimension, Algorithm~\ref{alg:squaring_algorithm}
may require significant computational resources due to the cubic (or quartic for ALS)
dependence on the separation rank $r$. 

Finally, we note that the applicability of our approach to global
optimization extends beyond the use of separated representations of
multivariate functions or tensors in CTD format. In fact, any structural
representation of a multivariate function that allows rapid computation
of the Hadamard product of corresponding tensors and has an associated
algorithm for reducing the complexity of the representation should
work in a manner similar to Algorithm~\ref{alg:squaring_algorithm}.

\bibliographystyle{plain}
\bibliography{common}

\end{document}